\magnification=\magstep1
\parskip 5pt
\font\Bbb=msbm10
\font\Bbbs=msbm8
\textfont12=\Bbb
\scriptfont12=\Bbbs
\font\rmsmall= cmr8
\font\bfone= cmbx10 scaled \magstep1
\let\mcd=\mathchardef
\mcd\Ee="7C45
\mcd\Pe="7C50
\mcd\Ze="7C5A
\def\al{\alpha}
\def\be{\beta}
\def\del{\delta}

\def\ga{\gamma}
\def\Ga{\Gamma}
\def\la{\lambda}
\def\bra{\langle}
\def\ket{\rangle}
\def\O{{\cal O}}
\def\nea{\nearrow}
\def\sea{\searrow}
\centerline{\bfone Random walk weakly attracted to a wall}
\bigskip
\centerline{Jo\"el De Coninck\footnote{$^{(1)}$}
{\rmsmall Centre de Recherche en Mod\'elisation Mol\'eculaire, Universit\'e de 
Mons-Hainaut, 20 Place du Parc, 7000 Mons, Belgium. 
Email: Joel.De.Coninck@crmm.umh.ac.be}, 
Fran\c cois Dunlop\footnote{$^{(2)}$}
{\rmsmall Laboratoire de Physique Th\'eorique et Mod\'elisation (CNRS - UMR 
8089), Universit\'e de Cergy-Pontoise, 95302 Cergy-Pontoise, France. 
Email: Francois.Dunlop@u-cergy.fr, Thierry.Huillet@u-cergy.fr}, 
Thierry Huillet$^{(2)}$}
\bigskip\noindent{\bf Abstract:} {\rm
We consider a random walk $X_n$ in $\Ze_+$, starting at $X_0=x\ge0$, with
transition probabilities
$$\Pe(X_{n+1}=X_n\pm1|X_n=y\ge1)={1\over2}\mp{\del\over4y+2\del}$$
and $X_{n+1}=1$ whenever $X_n=0$. We prove 
$\Ee X_n\sim{\rm const.}\,n^{1-{\del\over2}}$ as $n\nea\infty$ when 
$\del\in(1,2)$. The proof is based upon the Karlin-McGregor spectral
representation, which is made explicit for this random walk.}
 
\medskip\noindent {\rmsmall KEYWORDS: Random walk, orthogonal
polynomials, pinning, wetting} 
 
\bigskip\noindent
{\bf 1. Introduction}
\medskip\noindent
Random walks have been used in many different fields of physics,
economics, biology... Usually, it evolves in a translation invariant
environment or a random environment whose average is translation invariant.
Random walks in inhomogeneous environment, in particular with a reflecting or
attracting or repelling wall, have been used
to mimic the behaviour of a liquid interface on top of a solid
substrate, {\sl wetting} phenomena, and a variety of other phenomena
associated with the {\sl pinning} of an interface. 
This is the motivation for the model presented below, whose specific
form was chosen so as to derive rigorously and explicitly the result
of competition 
between a relatively long ranged attraction and reflection at the wall.

We consider a random walk $X_n$ in $\Ze_+$, defined by $X_0=x\ge0$,
$$\eqalign{
p_y=\Pe(X_{n+1}=X_n+1|X_n=y\ge1)&={1\over2}-{\del\over4y+2\del}
={y\over 2y+\del}\cr
q_y=\Pe(X_{n+1}=X_n-1|X_n=y\ge1)&={1\over2}+{\del\over4y+2\del}
={y+\del\over 2y+\del}} \eqno{(1.1)}
$$
and $X_{n+1}=1$ whenever $X_n=0$, i.e. $p_0=1$, $q_0=0$. Let
$\Pe_x(\cdot)$ denote the corresponding probabilities, and $\Ee_x(\cdot)$
the corresponding expectation values. 

The walk obeys the detailed balance condition with respect to the
measure $\pi$ on $\Ze_+$ defined up to a multiplicative constant by
$$
\pi_y=\pi_0\prod_{z=0}^{y-1}{p_z\over q_{z+1}}
=\pi_0{(2y+\del)\Ga(\del+1)\Ga(y)\over\Ga(y+\del+1)}\,,\quad y\ge1, 
\eqno{(1.2)}
$$
which obeys $\pi_y\sim y^{-\del}$ when $y\to\infty$. We restrict
our attention to $\del>1$, and normalise $\pi$ as a probability
measure, with $\pi_0=(\del-1)/(2\del)$.
The dynamics is periodic, with the odd and even components exchanged
under one step of the dynamics, $\pi=(\pi_{\rm even}+\pi_{\rm odd})/2$, with 
$\pi_{\rm even}=(2\pi_0,0,2\pi_2,\dots)$, and
$\pi_{\rm odd}=(0,2\pi_1,0,2\pi_3,\dots)$. 
Starting from $X_0=0$,
we have convergence in law $X_{2n}\to X_\infty^{\rm even}$ and
$X_{2n+1}\to X_\infty^{\rm odd}$. The first moment exists only for $\del>2$,
with 
$\Ee X_\infty=\Ee(X_\infty^{\rm even}+X_\infty^{\rm odd})/2
=\del/(2(\del-2))$.

We focus our attention to $1<\del<2$, so that $\Ee_x X_n\to\infty$
as $n\nea\infty$. For $|z|<1$ let
$$
g_e(z)=\sum_{n\,{\rm even}}z^n\Ee_0 X_n\,,\qquad
g_o(z)=\sum_{n\,{\rm odd}}z^n\Ee_0 X_n\eqno{(1.3)}
$$
Our main result is the following theorem:
\smallskip\noindent
{\bf Theorem 1}: {\sl Let $\del\in(1,2)$. Then as $z\nea1$ or $n\nea\infty$, 
$$
g_e(z)=\Ga\Bigl(2-{\del\over2}\Bigr)K_\del\,(1-z)^{{\del\over2}-2}
\biggl(1+\O\Bigl((1-z)^{(1-{\del\over2})(\del-1)}\Bigr)\biggr)\eqno{(1.4)}
$$
$$
\Ee_0X_n=K_\del\, n^{1-{\del\over2}}
\biggl(1+\O\Bigl(n^{-{(1-{\del\over2})(\del-1)}}\Bigr)\biggr)\eqno{(1.5)}
$$
where  $K_\del$ is a constant depending upon $\del$ as defined in (4.12).
}\smallskip
\noindent{\bf Remark 1}: The constant $K_\del\nea+\infty$ when
$\del\sea1$ or $\del\nea2$, but the $\O(\cdot)$ are not uniform in $\del$,
the theorem says nothing about $\del=1$ or $\del=2$.

\noindent{\bf Remark 2}: The odd generating function $g_o(z)$ obeys
the same bounds as $g_e(z)$. For any starting point $x$, $\Ee_xX_n$ obeys
the same bounds as $\Ee_0X_n$.

\noindent{\bf Remark 3}: (1.5) follows from (1.4) by a Tauberian
theorem [F, Thm. 5 p 447]. This uses monotonicity of $\Ee_0X_n$ for
$n$ even, which holds for any reflected random walk, as recalled below together
with more specific monotonicity arguments.

Consider two walks: $X_n$ started at $x$, with parameter $\del$, and
$X_n'$ started at $x'$, with parameter $\del'$. Assume $\del\ge\del'$ and
$x\le x'$ and $x'-x$ even. The two walks may be coupled, e.g. using a
single random number for both walks when they meet.
This implies that $\Ee_xX_n$ 
is for each $n$ an increasing function of $x$ over the even integers,
and also over the odd integers, and a decreasing function of
$\del\in(-1,+\infty)$. Monotonicity in $\delta$ is an example of a
general monotonicity property in the transition probabilities $p_y$, $q_y$. 
Then
$$
\Ee_0X_{n+2}=\Pe_0(X_2=0)\Ee_0X_n+\Pe_0(X_2=2)\Ee_2X_n\ge\Ee_0X_n\eqno{(1.6)}
$$
Therefore $\Ee_0X_n$ is an increasing function of $n$ over the even integers,
and also over the odd integers. And the same property holds for $\Ee_1X_n$.
But $\Ee_xX_n$ for $x\ge2$ is not monotonous in $n$, even or odd.
We have
$$
\Ee(X_{n+1}-X_n|X_n=y)=-{\del\over2y+\del}\,,\qquad y\ge1\eqno{(1.7)}
$$  
so that $\Ee_0X_{2n}<\Ee_0X_{2n-1}$. And since $\Ee_0X_n$ is
increasing over the even (or the odd) integers, we have
$\Ee_0X_{2n+1}>\Ee_0X_{2n}$. 

(1.4) will be proven in Section 4. In Section 2 we exhibit the
orthogonality measure of our random walk polynomials: this is Theorem
2, our main technical result, opening the way for Theorem 1. In Section 3 we
solve a differential equation for generating functions, differential equation
associated with the recursion formula for random walk polynomials.

Our random walk is of a very special form, but it could be used for
comparison or as input for more realistic models. For example the
polymer pinning model of Alexander and Zygouras [AZ] combines
i.i.d. disorder with a spatially inhomogeneous Markov chain, which
could be built from our random walk.

Other aspects of this model and closely related models will be treated
in a forthcoming paper [H].

\bigskip\noindent
{\bf 2. Random walk polynomials and their orthogonality measure}
\medskip\noindent

Our first tool is the Karlin-McGregor representation theorem [KM]:
let $L^2(\pi)$ denote the Hilbert space of complex sequences
$(f_y)_{y\in\Ze_+}$ obeying $\sum_{y=0}^\infty|f_y|^2\pi_y<\infty$. Then
$$
(Tf)_y=p_yf_{y+1}+q_yf_{y-1} \eqno{(2.1)}
$$
defines in $L^2(\pi)$ a self-adjoint operator $T$ of norm less or equal to one.
Let $e_0=(1,0,0,\dots)$. The Karlin-McGregor representation theorem gives 
$$
\Pe_x(X_n=y)={\pi_y\over\pi_0}\bra T^nQ_x(T)Q_y(T)e_0,e_0\ket \eqno{(2.2)}
$$
where $\{Q_y(t)\}_{y\in\Ze_+}$ is a family of polynomials of degree
$y$ in $t$, defined recursively by $Q_{-1}=0$, $Q_0=1$, and
$$
tQ_y=p_yQ_{y+1}+q_yQ_{y-1}\,,\quad y\ge0 \eqno{(2.3)}
$$
giving polynomials $Q_y$ of degree $y$ and parity $(-1)^y$, with
$Q_y(1)=1$ for all $y\ge0$. 
Using the spectral resolution $\{E_t\}$ of the self-adjoint operator $T$,
and $d\mu(t)=d\bra E_te_0,e_0\ket$, one gets
$$
\Pe_x(X_n=y)={\pi_y\over\pi_0}\int_{-1}^1d\mu(t)t^nQ_x(t)Q_y(t) \eqno{(2.4)}
$$
which implies that $\{Q_y(t)\}_{y\in\Ze_+}$ is a family of orthogonal
polynomials in the probability measure $d\mu(t)$,  
$$
\int_{-1}^1d\mu(t)Q_x(t)Q_y(t)={\pi_0\over\pi_y}\del_{x,y}\,,\quad x,y\ge0
 \eqno{(2.5)}
$$
Given the family $\{Q_y(t)\}_{y\in\Ze_+}$, (2.5) characterizes
a unique probability measure $d\mu$, termed the orthogonality measure
of the family. Letting $n\nea\infty$ in (2.4) shows that [KM, pp 70-71]
$$
d\mu(t)=\pi_0\Bigl(\del(t-1)+\del(t+1)\Bigr)\ +\ d\mu^c(t) \eqno{(2.6)}
$$
where $d\mu^c(t)$ is absolutely continuous with respect to the
Lebesgue measure. Indeed as $n\nea\infty$ with $n+y-x$ even, the LHS
of (2.4) tends to $2\pi_y$, while the contribution from $d\mu^c$ to the
RHS tends to zero.

Our second tool is Dette's theorem [D]:
the orthogonality measure $d\mu^1$ of the first associated polynomials of a
random walk is related to the orthogonality measure $d\mu^*$ of the dual random
walk through
$$
d\mu^1(t)={1\over q_1}(1-t^2)d\mu^*(t) \eqno{(2.7)}
$$
The first associated polynomials $Q_y^1$ are defined by $Q_{-1}^1=0$,
$Q_0^1=1$ and
$$
tQ_y^1=p_{y+1}Q_{y+1}^1+q_{y+1}Q_{y-1}^1\,,\quad y\ge0 \eqno{(2.8)}
$$
The dual random walk polynomials $Q_y^*$ are defined by $p_y^*=q_y$ and
$q_y^*=p_y$ for all $y$ except $p_0^*=1$ and $q_0^*=0$, so that
$Q_0^*=1$, $Q_1^*=t$ and
$$
tQ_y^*=p_y^*Q_{y+1}^*+q_y^*Q_{y-1}^*=q_yQ_{y+1}^*+p_yQ_{y-1}^*\,,\quad y\ge1 
\eqno{(2.9)} 
$$
Dette's theorem may be applied starting from the dual, giving
$$
d\mu^{*,1}(t)={1\over p_1}(1-t^2)d\mu(t) \eqno{(2.10)}
$$
The first associated dual polynomials $Q_y^{*,1}$ are defined by
$Q_{-1}^{*,1}=0$, $Q_0^{*,1}=1$ and
$$
tQ_y^{*,1}=q_{y+1}Q_{y+1}^{*,1}+p_{y+1}Q_{y-1}^{*,1}\,,\quad y\ge0 
\eqno{(2.11)}
$$
The general definitions (2.3)(2.8)(2.9)(2.11) are made explicit by inserting
our $p_y$'s and $q_y$'s, giving
$$\eqalign{
(2y+\del)tQ_y&=yQ_{y+1}+(y+\del)Q_{y-1}\,,\quad y\ge1\,,\quad Q_1(t)=t\cr
(2y+2+\del)tQ_y^1&=(y+1)Q_{y+1}^1+(y+1+\del)Q_{y-1}^1\,,\quad y\ge0 \cr
(2y+\del)tQ_y^*&=(y+\del)Q_{y+1}^*+yQ_{y-1}^*\,,\quad y\ge1 \cr
(2y+2+\del)tQ_y^{*,1}&=(y+1+\del)Q_{y+1}^{*,1}+(y+1)Q_{y-1}^{*,1}\,,\quad y\ge0
}\eqno{(2.12)} 
$$
We thus have four distinct families of orthogonal polynomials, and the
corresponding four distinct orthogonality measures $\mu$, $\mu^1$,
$\mu^*$, $\mu^{*,1}$.

The third and last step is to relate our polynomials to Gegenbauer polynomials
of index $\la$, defined  by $G_{-1}^\la=0$, $G_0^\la=1$ and
$$
(2y+2\la)tG_y^\la=(y+1)G_{y+1}^\la+(y-1+2\la)G_{y-1}^\la\,,\quad y\ge0
\eqno{(2.13)} 
$$
It appears that $Q_y^1=G_y^{{\del\over2}+1}$, so that
$Q_y^{1,*}=G_y^{{\del\over2}+1,*}$, and also
$Q_y^*={y!\Ga(\del)\over\Ga(y+\del)}G_y^{\del\over2}$. These in turn satisfy
$$
G_y^{{\del\over2}+1,*}(t)={(y+1)!\Gamma(\del+1)\over\Gamma(y+1+\del)}
G_y^{{{\del\over2},1}}(t)
={(y+1)!\Gamma({\del+3\over2})\over\Gamma(y+{\del+3\over2})}
P_y^{{\del-1\over2},{\del-1\over2}}(t;1) \eqno{(2.14)}
$$
where $P_x^{\al,\be}(t;c)$ are $c$-associated Jacobi polynomials,
or Wimp polynomials [W].
The orthogonality measure $d\mu^{*,1}$ is therefore also the
orthogonality measure of the $P_y^{{\del-1\over2},{\del-1\over2}}(t;1)$
polynomials, namely [W, Th. 3 p. 996]
$$
d\mu^{*,1}(t)={(1-t^2)^{\del-1\over2}\over|F(t)|^2}dt\Big/{\rm normalisation}
\eqno{(2.15)}
$$
where
$$\eqalign{
F(t)&= {_2F_1}\Bigl(1,1-\del;{3-\del\over2};{1+t\over2}\Bigr)
+Ke^{i\pi{\del-1\over2}}\Bigl({1+t\over2}\Bigr)^{\del-1\over2}\;
{_2F_1}\Bigl({1+\del\over2},{1-\del\over2};{1+\del\over2};{1+t\over2}\Bigr)\cr 
&= {_2F_1}\Bigl(1,1-\del;{3-\del\over2};{1+t\over2}\Bigr)
+Ke^{i\pi{\del-1\over2}}\,\Bigl({1-t^2\over4}\Bigr)^{\del-1\over2}\cr 
&=-1+\O\bigl((1-t)^{\del-1\over2}\bigr)
\qquad{\rm as}\quad t\to1 }\eqno{(2.16)}
$$
with $_2F_1$ the Gauss hypergeometric function and
$$
K={\Ga(\del)\Ga({1-\del\over2})\over\Ga({\del-1\over2})}\eqno{(2.17)}
$$
We have $F(-t)=-F(t)^*$. Then (2.6)(2.10) and (2.15) give
$$
d\mu^c(t)={(1-t^2)^{\del-3\over2}\over|F(t)|^2}dt\Big/{\rm normalisation}
\eqno{(2.18)}
$$
which yields:
\smallskip\noindent{\bf Theorem 2}: {\sl The orthogonality measure of the
family of polynomials defined by (2.3) with (1.1) and $1<\del<2$ is the even
probability measure on $[-1,1]$ defined by
$$
d\mu(t)=\pi_0\Bigl(\del(t-1)+\del(t+1)\Bigr)\ +\ d\mu^c(t) \eqno{(2.19)}
$$
where $\del(\cdot)$ is the Dirac measure at 0 and $d\mu^c$ is given by
(2-16)-(2.18) and $\pi_0={\del-1\over2\del}$, with the normalisation
$\int_{-1}^1d\mu^c=\del^{-1}$.} 

\vfill\eject
\parskip 2pt
\noindent{\bf 3. Generating function}
\medskip\noindent
Using (2.4) and $X_n\le X_0+n$, we have
$$
\Ee_0 X_n=\sum_{y=1}^n{y\pi_y\over\pi_0}\int_{-1}^1d\mu(t)t^nQ_y(t)\eqno{(3.1)}
$$ 
For $n$ even, using $Q_y$ orthogonal to $Q_0\equiv1$, and for $y$ odd
$Q_y$ also orthogonal to $t^n$, and then using (2.6), we have  
$$\eqalign{
\Ee_0 X_n
&=\sum_{y=2\atop{\rm even}}^n{y\pi_y\over\pi_0}
\int_{-1}^1d\mu^c(t)(t^n-1)Q_y(t)\cr
&=2\int_0^1d\mu^c(t)(t^n-1)\sum_{y=0\atop{\rm even}}^n{y\pi_y\over\pi_0}Q_y(t)
}\eqno{(3.2)}
$$
The corresponding generating function is defined as 
$$\eqalign{
g_e(z)&=\sum_{n\,{\rm even}}z^n\Ee_0 X_n\cr
&=2\int_0^1d\mu^c(t)\sum_{n\,{\rm even}}z^n(t^n-1)
\sum_{y=0\atop{\rm even}}^n{y\pi_y\over\pi_0}Q_y(t)\cr
&=2\int_0^1d\mu^c(t)\sum_{y\,{\rm even}}{y\pi_y\over\pi_0}Q_y(t)
\Bigl({(zt)^y\over1-z^2t^2}-{z^y\over1-z^2}\Bigr)}\eqno{(3.3)}
$$ 
For $n$ odd, using $Q_y$ orthogonal to $Q_1\equiv t$ for $y\ge2$, and
(2.5) for $y=1$,
$$\eqalign{
\Ee_0 X_n&=\sum_{y=1}^n{y\pi_y\over\pi_0}\int_{-1}^1d\mu(t)(t^n-t)Q_y(t)
+{\pi_1\over\pi_0}\int_{-1}^1d\mu(t)tQ_1(t)\cr
&=1+2\int_0^1d\mu^c(t)(t^n-t)
\sum_{y=1\atop{\rm odd}}^n{y\pi_y\over\pi_0}Q_y(t)
=\Ee_1X_{n-1}}\eqno{(3.4)}
$$ 
and the corresponding generating function
$$\eqalign{
g_o(z)&=\sum_{n\,{\rm odd}}z^n\Ee_0 X_n\cr
&={z\over1-z^2}+2\int_0^1d\mu^c(t)\sum_{n\,{\rm odd}}z^n(t^n-t)
\sum_{y=1\atop{\rm odd}}^n{y\pi_y\over\pi_0}Q_y(t)\cr
&={z\over1-z^2}+2\int_0^1d\mu^c(t)\sum_{y\,{\rm odd}}{y\pi_y\over\pi_0}Q_y(t)
\Bigl({(zt)^y\over1-z^2t^2}-{tz^y\over1-z^2}\Bigr)}\eqno{(3.5)}
$$ 
The recursion (2.3)(2.12) defining the random walk polynomials gives
a poor uniform bound for these polynomials, e.g. $|Q_y(t)|<3^y$ for
$t\in(-1,1)$. We define, for $|u|<1/3$,
$$
\Psi_t(u)=\sum_{y=1}^\infty{\pi_y\over\pi_0}Q_y(t)u^y
=\Ga(\del+1)\sum_{y=1}^\infty(2y+\del){\Ga(y)\over\Ga(y+\del+1)}Q_y(t)u^y
\eqno{(3.6)}
$$
and aim at an analytic continuation giving, for $|z|<1$,
$$\eqalign{
g_e(z)=z\int_0^1d\mu^c(t)\biggl[\biggl({t\Psi_t'(zt)\over1-z^2t^2}
-{\Psi_t'(z)\over1-z^2}\biggr)
-\biggl({t\Psi_t'(-zt)\over1-z^2t^2}-{\Psi_t'(-z)\over1-z^2}\biggr)\biggr]\cr
g_o(z)={z\over1-z^2}
+z\int_0^1d\mu^c(t)t\biggl[\biggl({\Psi_t'(zt)\over1-z^2t^2}
-{\Psi_t'(z)\over1-z^2}\biggr)
+\biggl({\Psi_t'(-zt)\over1-z^2t^2}-{\Psi_t'(-z)\over1-z^2}\biggr)\biggr]
}\eqno{(3.7)}
$$
The recursion (2.3)(2.12) may be converted into a first order differential
equation for $\Psi_t(u)$. We first get a differential equation for
$$
\Phi_t(u)=\Ga(\del+1)\sum_{y=1}^\infty{\Ga(y)\over\Ga(y+\del+1)}Q_y(t)u^y
\equiv\sum_{y=1}^\infty H_y(t)u^y\eqno{(3.8)}
$$
with
$$
H_y(t)={1\over2y+\del}{\pi_y\over\pi_0}Q_y(t)
={\Ga(\del+1)\Ga(y)\over\Ga(y+\del+1)}Q_y(t)\,,\qquad y\ge1\eqno{(3.9)}
$$
obeying, with any arbitrary value for $H_0$,
$$
(2y+\del)tH_y=(y+\del+1)H_{y+1}+(y-1)H_{y-1}\,,\qquad y\ge2\eqno{(3.10)}
$$
with
$$
H_1(t)={t\over\del+1}\,,\qquad H_2(t)={Q_2(t)\over(\del+1)(\del+2)}
={(\del+2)t^2-(\del+1)\over(\del+1)(\del+2)}\eqno{(3.11)}
$$
Multiplying (3.10) by $u^y$ and summing over $y\ge2$ yields
$$
(1-2tu+u^2)\Phi'_t(u)=t-u-\del(u^{-1}-t)\Phi_t(u)\eqno{(3.12)}
$$
with $\Phi_t(0)=0$, whose solution is
$$
\Phi_t(u)=\del^{-1}-u^{-\del}(1-2tu+u^2)^{\del\over2}
\int_0^udv\,v^{\del-1}(1-2tv+v^2)^{-{\del\over2}}\eqno{(3.13)}
$$
We thus get:
\smallskip\noindent{\bf Lemma 3}: {\sl The function $\psi_t(u)$ defined in
(3.6) may be expressed as 
$$
\Psi_t(u)=2u\Phi'_t(u)+\del\Phi_t(u)
=-1+{\del(1-u^2)\over1-2tu+u^2}\Bigl(\del^{-1}-\Phi_t(u)\Bigr)\eqno{(3.14)}
$$
where $\Phi_t(u)$ is the solution of the differential equation (3.12). It
extends to an analytic function in the disc $|u|<1$. Its derivative
may be expressed as
$$
\Psi_t'(u)={\del(1-u^2)\over u\,(1-2tu+u^2)}
-{B_t(u)\,\del\over u^\del\,(1-2tu+u^2)^{2-{\del\over2}}}
\int_0^u{dv\,v^{\del-1}\over(1-2tv+v^2)^{\del\over2}}\eqno{(3.15)}
$$
with
$$
B_t(u)=4u(1-t)-2t(1-u)^2+{\del(1-ut)(1-u^2)\over u}\eqno{(3.16)}
$$
}

\bigskip\noindent
{\bf 4. Proof of (1.4)}
\medskip\noindent
The leading order and next to leading order in (3.7) as $z\nea1$ are found in
$$
g_1(z)=\int_0^1d\mu_c(t)\biggl({zt\Psi_t'(zt)\over1-z^2t^2}
-{z\Psi_t'(z)\over1-z^2}\biggr)\eqno{(4.1)}
$$
Indeed the leading orders come from the singularity at $t=1$ in the integral.
In (3.7) there is a symmetry or anti-symmetry as $z\to-z$ and $t\to-t$ jointly,
associated with the even/odd symmetries. The terms not included in (4.1)
correspond to $-z\simeq-1$, not singular with $t>0$.
For $1-z\ll1$ and $1-t\ll1$ we have
$$\eqalign{
{t\over1-z^2t^2}-{1\over1-z^2}&=-(1-t){1+z^2t^2\over(1-z^2t^2)(1-z^2)}
\sim-{1\over2}{1-t\over1-z}{1\over(1-z)+(1-t)}\cr
{t\over1-z^2t^2}+{1\over1-z^2}&={t(1-z^2)+1-z^2t^2\over(1-z^2t^2)(1-z^2)}
\sim{1\over2}{1\over1-z}{2(1-z)+1-t\over(1-z)+(1-t)}
}\eqno{(4.2)}
$$
(3.15)(3.16) may be written as
$$
q_t(z)^2\Psi_t'(z)=\del{1-z^2\over z}q_t(z)
-B_t(z)\del z^{-\del}q_t(z)^{\del\over2}
\int_0^zdv\,v^{\del-1}q_t(v)^{-{\del\over2}}\eqno{(4.3)}
$$
with
$$
q_t(z)=1-2tz+z^2=(1-z)^2+2z(1-t)\,,\qquad q_t'(z)=2(z-t)\eqno{(4.4)}
$$
Together with (4.2), we have to estimate $\Psi_t'(zt)\pm\Psi_t'(z)$
as $z\nea1$, $t\nea1$. We have
$$
{1-z^2t^2\over zt}\sim{1-z^2\over z}\Bigl(1+{1-t\over1-z}\Bigr)\eqno{(4.5)}
$$
$$\eqalign{
B_t(z)&\sim4(1-t)-2(2-\del)(1-z)(1-t)+2(\del-1)(1-z)^2\cr
B_t'(z)&\sim2(2-\del)(1-t)-4(\del-1)(1-z)\cr
B_t(zt)&\sim B_t(z)\Bigl(1-(1-t){B_t'(z)\over B_t(z)}\Bigr)
}\eqno{(4.6)}
$$
$$\eqalign{
\int_0^1{dv\,v^{\del-1}\over q_t(v)^{\del\over2}}
&=\int_0^1{dv\over\bigl((1-v)^2+2(1-t)\bigr)^{\del\over2}}
+\int_0^1dv\biggl[\,{v^{\del-1}\over q_t(v)^{\del\over2}}
-{1\over\bigl((1-v)^2+2(1-t)\bigr)^{\del\over2}}\biggr]\cr
&=(1-t)^{1-\del\over2}\int_0^{(1-t)^{-{1\over2}}}
{dx\,\over\bigl(x^2+2\bigr)^{\del\over2}}\ +\ \O(1)\cr
&=(1-t)^{1-\del\over2}\,
{2^{-{1+\del\over2}}\sqrt\pi\Ga({\del-1\over2})\over\Ga({\del\over2})}
\ +\ \O(1)
}\eqno{(4.7)}
$$
The range $(0,1)$ of the integral in (4.1) is split into four
intervals, according to
$$
0<1-(1-z)^\al<1-(1-z)^\be<1-(1-z)^\ga<1\eqno{(4.8)}
$$
with $0<\al<1<\be<2<\ga$ and respective contributions denoted 
$g_1^{\al}\,,g_1^{\al\be}\,, g_1^{\be\ga}\,,  g_1^{\ga}$.
We begin with the leading contribution, $g_1^{\al\be}$:
\item{$\bullet$} $1-(1-z)^{\al}<t<1-(1-z)^{\be}\,$:

\noindent 
$q_t(z)\sim 2(1-t)$ and $B_t(z)\sim4(1-t)$. Then (4.3) with (4.7) yields 
$$
(1-t)^2\Psi_t'(z)\sim\del(1-z)(1-t)-(1-t)^{3\over2}\del\sqrt{\pi\over2}
{\Ga({\del-1\over2})\over\Ga({\del\over2})}\eqno{(4.9)}
$$
and
$$\eqalign{
\int_{1-(1-z)^{\al}}^{1-(1-z)^\be}dt\,(1-t^2)^{\del-3\over2}
\Bigl({t\over1-z^2t^2}-{1\over1-z^2}\Bigr)&\Psi_t'(z)\sim\cr
\sim2^{\del-5\over2}\del\sqrt{\pi\over2}
{\Ga({\del-1\over2})\over\Ga({\del\over2})}
\int_{1-(1-z)^{\al}}^{1-(1-z)^\be}&dt\,
(1-t)^{{\del\over2}-2}\,{1-t\over1-z}{1\over(1-z)+(1-t)}\cr
=2^{\del-5\over2}\del\sqrt{\pi\over2}
{\Ga({\del-1\over2})\over\Ga({\del\over2})}
(1-z)^{{\del\over2}-2}&\int_{(1-z)^{\be-1}}^{(1-z)^{\al-1}}dx\,
{x^{{\del\over2}-1}\over1+x}\cr
=2^{\del-5\over2}\del\sqrt{\pi\over2}\Ga\Bigl({\del-1\over2}\Bigr)
\Ga\Bigl(1-{\del\over2}\Bigr)(1-z)^{{\del\over2}-2}
\Bigl(&1+\O\Bigl((1-z)^{(1-\al)(1-{\del\over2})}
+(1-z)^{(\be-1){\del\over2}}\Bigr)\Bigr)
}\eqno{(4.10)}
$$
Collecting all previous error terms we get
$$\eqalign{
g_1^{\al\be}(z)=\Ga\Bigl(2-{\del\over2}\Bigr)K_\del(1-z)^{{\del\over2}-2}
\Bigl[1&+\O\Bigl((1-z)^{(1-\al)(1-{\del\over2})}
+(1-z)^{\al{\del-1\over2}}+(1-z)^{\al\over2}\cr
&+(1-z)^{(\be-1){\del\over2}}+(1-z)^{2-\be}\Bigr)\Bigr]
}\eqno{(4.11)}
$$
with
$$
K_\del=2^{\del-5\over2}\sqrt{\pi\over2}{\Ga\Bigl({\del-1\over2}\Bigr)
\over1-{\del\over2}}
\bigg/\int_{-1}^1{(1-t^2)^{\del-3\over2}\over|F(t)|^2}dt\eqno{(4.12)}
$$
and $F(t)$ given by (2.16). The five error terms in $\O(\cdots)$ in
(4.11) come respectively from:  
\item{1)} Extending the $x$-integral to $+\infty$ in (4.10)
\item{2)} $F(t)=-1+\O\bigl((1-t)^{\del-1\over2}\bigr)$
\item{3)} $\Psi_t'(zt)-\Psi_t'(z)\sim\del$, derived from (4.3)-(4.6),
with the leading contribution from (4.5).
\item{4)} Extending the $x$-integral to $0$ in (4.10)
\item{5)} $q_t(z)\sim 2(1-t)\Bigl(1+{(1-z)^2\over2(1-t)}\Bigr)$ and
$B_t(z)\sim 4(1-t)\Bigl(1+(\del-1){(1-z)^2\over2(1-t)}\Bigr)$.

Choosing 
$$
\al=\sup_{\del\in(1,2)}\Bigl(\min\Bigl\{(1-\al)(1-{\del\over2}),
\al{\del-1\over2}  \Bigr\}\Bigr)=2-\del\eqno{(4.13)}
$$ 
and $\be$ within $4-2\sqrt2<\be<2$,
e.g. $\be=3/2$, then brings (4.11) to (1.4). The remaining ranges in
the $t$-integral also fall within the previous error estimates:
\item{$\bullet$} $0<t<1-(1-z)^{\al}\,$: 
$$
g_1^{\al}=(1-z)^{{\del\over2}-2}\,
\O\Bigl((1-z)^{(1-\al)(1-{\del\over2})}+(1-z)^{\al\over2}
+(1-z)^{\al{\del-1\over2}}\Bigr)\eqno{(4.14)}
$$
\item{$\bullet$} $1-(1-z)^{\be}<t<1-(1-z)^{\ga}\,$:
$$
g_1^{\be\ga}=(1-z)^{{\del\over2}-2}\,
\O\Bigl((1-z)^{(\be-1){\del\over2}}+(1-z)^{2-\be}\Bigr)\eqno{(4.15)}
$$
\item{$\bullet$} $1-(1-z)^{\ga}<t<1\,$:
$$
g_1^{\ga}=\O\Bigl((1-z)^{-1}\Bigr)\eqno{(4.16)}
$$
This completes the proof of (1.4) and Theorem 1.

\bigskip\noindent
{\bf Acknowledgments}: The authors thank Yann Costes (Service d'Informatique
Recherche, UCP) for help with {\sl Mathematica} for some checks and error
correcting work. F. D. and T. H. acknowledge kind hospitality at
Universit\'e de Mons-Hainaut and CRMM where the present work was initiated.
\bigskip\noindent
{\bf References}
\medskip\noindent

\noindent\item{[AZ]} K.S. Alexander, N. Zygouras: {\sl 
Quenched and Annealed Critical Points in Polymer Pinning Models},
http://arxiv.org/abs/0805.1708

\noindent\item{[D]}  H. Dette: {\sl First return probabilities of
birth and death chains and associated orthogonal polynomials},
Proc. Amer. Math. Soc. {\bf 129}, 1805--1815 (2001).

\noindent\item{[F]} W. Feller: {\sl An introduction to probability
theory and its applications},  Vol. II (Wiley, 1971). 

\noindent\item{[H]} Th. Huillet: {\sl 
Random walk with long-range interaction with a
barrier and its dual: Exact results}, in preparation. 

\noindent\item{[KM]}  S. Karlin, J. McGregor: {\sl  Random walks}, Illinois
J. Math. {\bf 3}, 66--81 (1959).

\noindent\item{[W]}  J. Wimp: {\sl  Explicit formulas for the associated Jacobi
polynomials and some applications}, Canad. J. Math. {\bf 39}, 983-1000 (1987).

\bye